\documentclass[11pt]{article}
\usepackage{fullpage}

\usepackage{amssymb,latexsym,amsmath,epsfig,amsthm}

\usepackage{mathtools}
\DeclarePairedDelimiter\abs{\lvert}{\rvert}%
\usepackage[loadshadowlibrary,textwidth=2cm,textsize=small,shadow]{todonotes}

\makeatletter
\let\oldabs\abs
\def\abs{\@ifstar{\oldabs}{\oldabs*}}

\usepackage{enumitem}

\usepackage{hyperref}
\hypersetup{
  pdflang=en-US,
  pdftitle={Visualizing the Sum-Product Conjecture},
  pdfauthor={Dr Kevin O'Bryant}
}
\newcommand{\seqnum}[1]{\href{https://oeis.org/#1}{\rm \underline{#1}}}

\DeclareMathOperator{\NSPP}{NSPP}
\DeclareMathOperator{\SPP}{SPP}
\DeclareMathOperator{\lcm}{LCM}

\newcommand{\NN}{\ensuremath{\mathbb N}}
\newcommand{\ZZ}{\ensuremath{\mathbb Z}}

\newcommand{\RR}{\ensuremath{\mathbb R}}

\newcommand{\floor}[1]{\lfloor {#1} \rfloor}

\newcommand{\ceiling}[1]{\lceil {#1} \rceil}

\newcommand{\totalsize}{\ensuremath{1\,162\,868}} 

\newtheorem{theorem}{Theorem}

\newtheorem{corollary}{Corollary}

\theoremstyle{definition}

\newtheorem{conjecture}{Conjecture}
\newtheorem*{spc}{Sum-Product Conjecture}

\title{Visualizing the Sum-Product Conjecture}
\author{
  Kevin O'Bryant\\
  Department of Mathematics, City University of New York,\\
  The College of Staten Island and The Graduate Center, NY, USA\\
  \href{mailto:kevin.obryant@csi.cuny.edu}{kevin.obryant@csi.cuny.edu}
  }
\date{December 7, 2024}

\begin{document}
\maketitle

\begin{abstract}
Let $\SPP(n)$ be the set $\left\{\big(|A+A|,|A A|\big) : A\subseteq {\mathbb N}, |A|=n\right\}$ of sum-product pairs, where $A+A$ is the sumset $\{a+b : a,b\in A\}$ and $A A$ is the product set $\{ab:a,b\in A\}$.
We construct a dataset consisting of \totalsize\ sets whose sum-product pairs are at least $84\%$ of $\SPP(n)$ for each $n\le 32$.
Notably, we do {\bf not} see evidence in favor of Erd\H{o}s's Sum-Product Conjecture in our dataset. For $n\le 6$, we prove the exact value of $\SPP(n)$. We include a number of conjectures, open problems, and observations motivated by this dataset, a large number of color visualizations.
\end{abstract}


\section{Introduction}
About eighteen months after September $1974$, P.~Erd\H{o}s~\cite{1976.Erdos} wrote:
\begin{quote}
  Let $a_1,\ldots,a_n$ an be a sequence of $n$ numbers. Consider the numbers
  \mbox{$a_i+a_j, a_i a_j$}. I conjectured about $18$ months ago that there are more than $n^{2-\epsilon}$ distinct numbers amongst them. The conjecture presumably holds whether the $a$'s are integers or real or complex numbers, but I have not succeeded in proving the much weaker inequality $n^{1+\epsilon}$, even when the $a$'s are integers.
\end{quote}
This has become known as the Erd\H{o}s Sum-Product Conjecture\footnote{It is more commonly described as a conjecture of Erd\H{o}s and Szemer\'edi, but this is not correct. In 1983, Erd\H{o}s and Szemer\'edi~\cite{ErdosSzemeredi} made the first substantial progress on the conjecture, proving the $n^{1+\epsilon}$ bound for positive integers with an unspecified positive $\epsilon$.}. A more formal statement of the conjecture follows. Let $A+A \coloneqq \{x+y : x,y\in A\}$ be the sumset and $AA \coloneqq \{xy : x,y\in A\}$ be the product set.
\begin{spc}
  For all $\epsilon>0$ there is an $n_0$ such that if $A$ is any set of positive integers with $|A|\ge n_0$, then
  \[ \big| (A + A)\cup(AA) \big| \ge |A|^{2-\epsilon}.\]
\end{spc}
\noindent Erd\H{o}s offered \$100 for a proof or disproof, and offered \$250 for a more exact bound~(see Problem F18~\cite{UPINT}).

After examining several trillion sets of positive integers, we have compiled a dataset of sum-product pairs
\[ \SPP(n) \coloneqq \left\{ \big(|A+A|,|AA|\big) : A \subseteq \NN, |A|=n\right\}\]
with $n\le 32$, and a set with the minimal maximum that we have found for each triple $|A|,|A+A|,|AA|$. While the dataset is certainly incomplete, we hope that it is a useful resource and is available at \url{http://www.math.csi.cuny.edu/obryant/SumProduct}. In this work, we present visualizations of the dataset and prove some of the basic features seen, and make several conjectures supported by the dataset. We prove exact values of $\SPP(n)$ for $n\le 6$.

Using the $\SPP(n)$ notation, we restate the Sum-Product Conjecture as follows.
\begin{spc}
  For every $\epsilon>0$, if $n>n_0$ then $\SPP(n)$ and $[0,n^{2-\epsilon}]^2$ are disjoint.
\end{spc}

Supported by our dataset, we make numerous conjectures, including the following two.
\begin{conjecture}\label{conj:SV}
  Let $A$ be a set of $n$ positive integers. Then
  \[|A+A||AA|^2 \ge \frac{n(n+1)}2\cdot (2n+1)^2 
  \]
  with equality if an only if $A$ is a geometric progression.
\end{conjecture}
\noindent There are sets of positive real numbers for which the inequality in Conjecture~\ref{conj:SV} does not hold.
\begin{conjecture}\label{conj:|A+A|+|AA|}
  Let $A$ be a set of $n\ge2$ positive integers. Then
  \[|A+A|+|AA| \ge |A|^{(1+\sqrt5)/2}.\]
\end{conjecture}
\begin{conjecture}\label{conj:unnormalized}
  For every $(x,y)\in[0,\frac12]^2$, there is a sequence of sets of positive integers $A$ with $|A|\to\infty$, $|A+A|\to x |A|^2$, and $|AA| \to y |AA|^2$.
\end{conjecture}
\noindent We discuss evidence for and against these conjectures at the appropriate places below.

\begin{figure}
  \includegraphics[width=\textwidth]{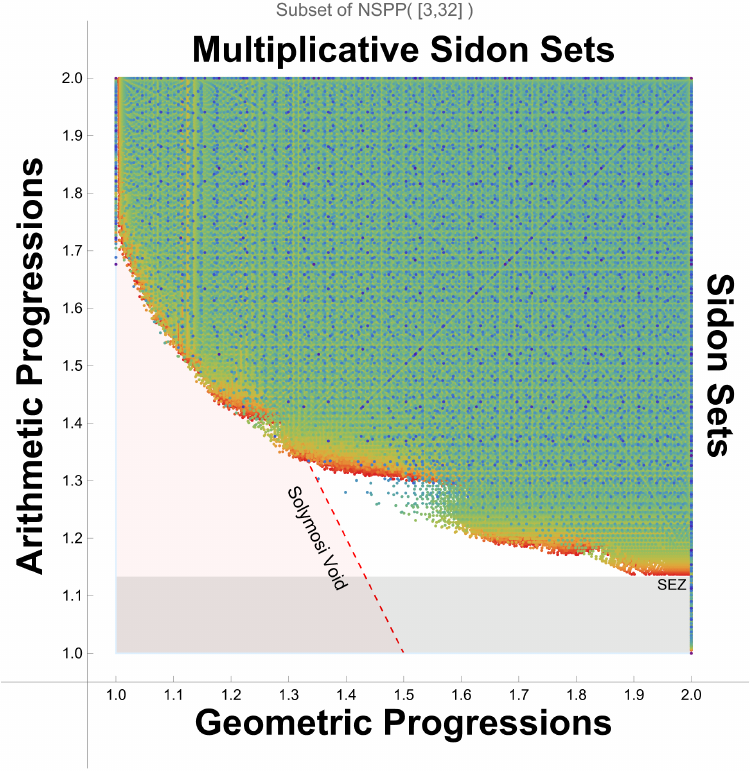}

  \caption{The set $\NSPP([3,32])$. Each point represents $(|A+A|,|AA|)$ for some set $A$, normalized logarithmically to lie in $[1,2]^2$. The sets range in size from $3$ to $32$, and the color of a point reflects its cardinality, with larger sets using the red end of the rainbow. The points reflecting sets with $n=32$ elements were laid down first, then $n=31$, and so on, so that purpler points from small $n$ typically conceal redder points beneath them.
  The Solymosi Void is a region known to not contain any limit points (as $n\to\infty$), and SEZ is a region known not to contain any (limit or otherwise) points of $\NSPP([3,32])$.}
  \label{fig:everythingK}
\end{figure}

\begin{figure}
  \includegraphics[width=\textwidth]{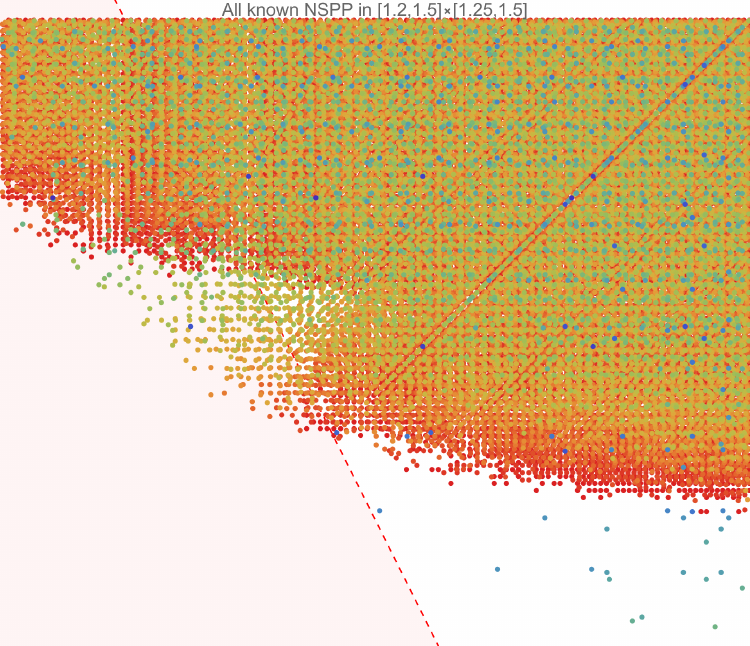}
  \caption{The set $\NSPP([3,32])$, restricted to $[1.2,1.5]\times[1.25,1.5]$. The color of the point indicates the size of the set with that normalized sum-product pair, with the reddest points corresponding to $n=32$. The $n=32$ points were laid down first, then $n=31$, and so on, so that points from small $n$ typically conceal redder points beneath them.}
  \label{fig:everythingzoom}
\end{figure}

\begin{figure}
  \begin{tabular}{cc}
    \includegraphics[width=0.47\textwidth]{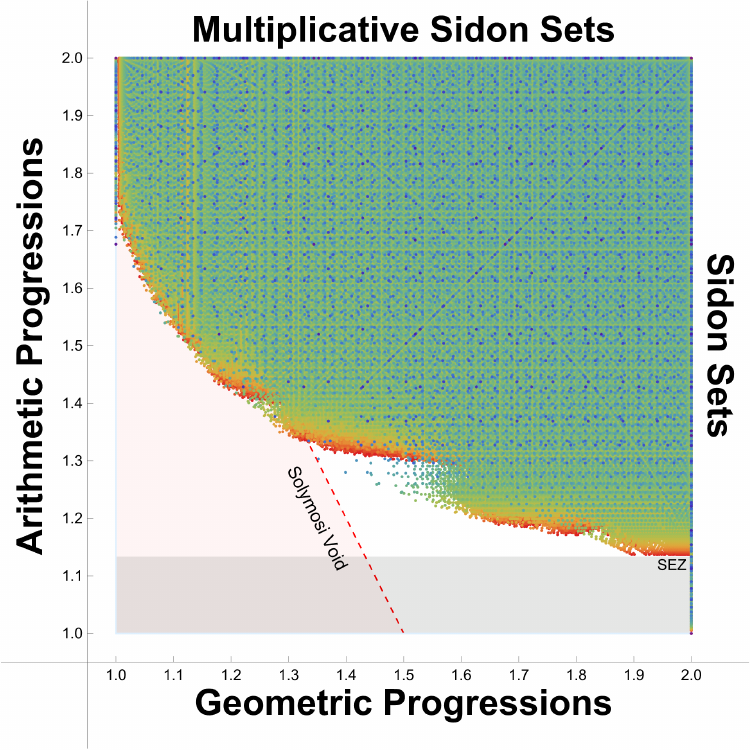} &  \includegraphics[width=0.47\textwidth]{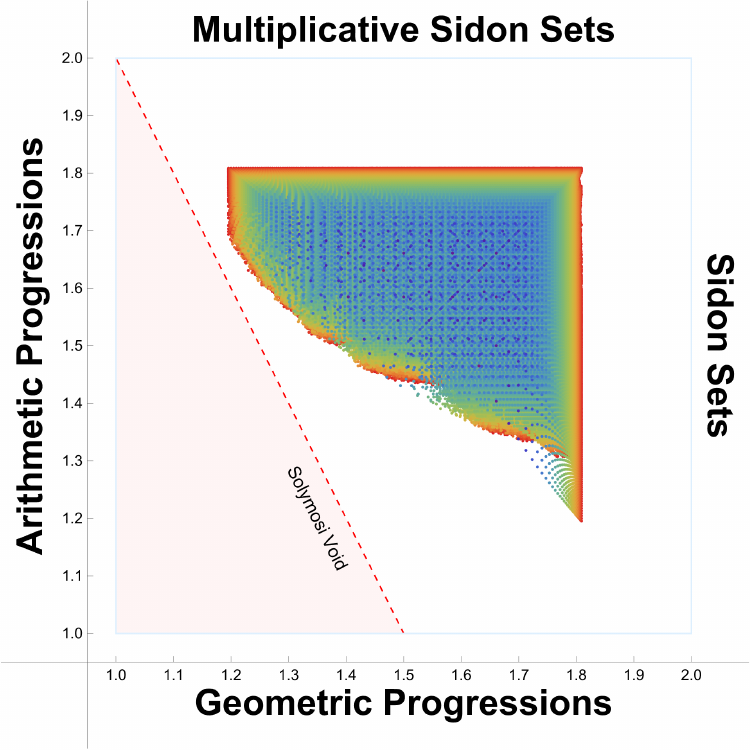} \\
    \includegraphics[width=0.47\textwidth]{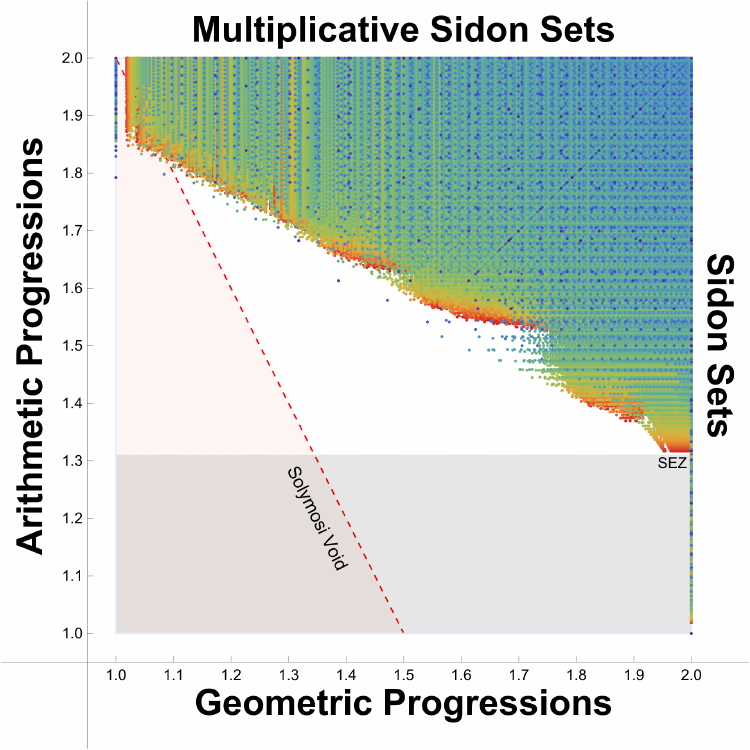} & \includegraphics[width=0.47\textwidth]{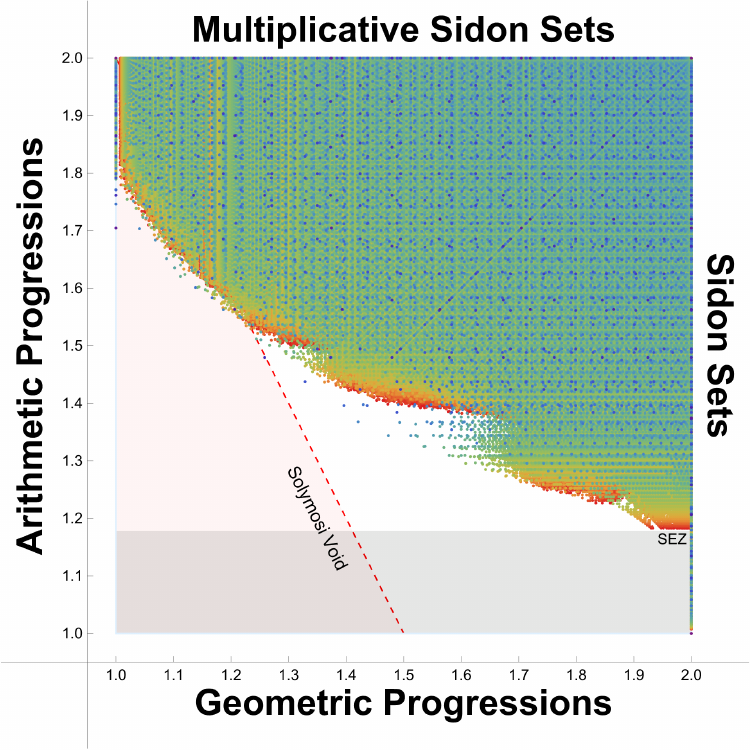}\\
  \end{tabular}
  \caption{The normalized sum-product pairs under normalizations $K,L,K^{(1)},K^{(2)}$. The reddish dots along the lower envelopes suggest that the Sum-Product Conjecture may be false, or at least the sum-product effect is not visible in sets with at most $32$ elements.}
  \label{fig:everything4}
\end{figure}

The most interesting images in this work are Figures~\ref{fig:everythingK}, \ref{fig:everythingzoom}, and \ref{fig:everything4}, which show the entire dataset but require some explanation. It is elementary that $2n-1\le |A+A|,|AA| \le n(n+1)/2$, so that $\SPP(n)$ is contained in $[2n-1,n(n+1)/2]^2$ and has cardinality at most $(n^2-3n+4)^2/4$. The Sum-Product Conjecture pulls our attention to $\log|A+A|,\log|AA|$, and in order to show $\SPP(n)$ for different $n$ on the same image, it is necessary to normalize the sum-product pairs. For $n\ge 3$, we define $4$ normalizing functions:
\begin{align*}
  K_n(x) &\coloneqq \log_n(x) +m_n x+b_n, & m_n,b_n\in\RR,\quad K_n(2n-1)=1,\quad K_n(n(n+1)/2)=2 \\
  L_n(x) &\coloneqq \log_n(x) & \\
  K_n^{(2)}(x) & \coloneqq \log_n(m_n x + b_n), & m_n,b_n\in\RR,\quad K_n^{(2)}(2n-1)=1, \quad K_n^{(2)}(n(n+1)/2)=2 \\
  K_n^{(3)}(x) & \coloneqq m_n \log_n(x) + b_n, & m_n,b_n\in\RR,\quad K_n^{(3)}(2n-1)=1, \quad K_n^{(3)}(n(n+1)/2)=2
\end{align*}
The author finds $K_n(x)$ the most natural, while $L_n(x)$ is more combinatorial. We include $K_n^{(2)}$ and $K_n^{(3)}$ to prevent mistaking an artifact of the normalization for a trait of the data.

In these terms, as $n\to\infty$ we have
  \[|A+A| = n^{K_n(|A+A|)+o(1)},\qquad |AA| = n^{K_n(|AA|)+o(1)},\]
and $K_n(|A+A|),K_n(|AA|) \in [1,2]^2$ for every set $A$ with $|A|\ge 3$. We define the sets (not multi-sets)
\[ \NSPP(j) \coloneqq \bigg\{ \big( K_n(\abs{A+A}), K_n(\abs{AA}) \big) : A \subseteq \NN, \abs{A}= j \bigg\} \]
and the multisets of normalized sum-product pairs
\[ \NSPP(J) \coloneqq \bigcup_{j\in J} \NSPP(j).\]
A more detailed explanation of the images in Figures~\ref{fig:everythingK} and~\ref{fig:everythingzoom} is in Subsection~\ref{subsec:everything pics}.

We restate the Sum-Product Conjecture using the $\NSPP$ notation.
\begin{spc}
  The limit points of the multiset $\NSPP([3,\infty))$ are confined to the two line segments $\{2\}\times [1,2]$, $[1,2]\times \{2\}$.
\end{spc}
\noindent
Loosely speaking, the Sum-Product Conjecture says that $\NSPP(j)$ moves north and east as $j\to \infty$.
As a corollary of a recent result of Nathanson~\cite{Nathansonsumsetsizes}, we know that every point in $\{2\}\times [1,2]$, $[1,2]\times \{2\}$ is a limit point of $\NSPP([3,\infty))$.

In Figure~\ref{fig:everythingK}, one sees $\NSPP([3,32])$ with points arising from larger $n$ colored more red. Points for the largest $n$ are laid down first, so that smaller $n$ (more purple) appear on top. Contrary to what the Sum-Product Conjecture suggests, as $|A|$ becomes larger we see points $\big( K_n(|A+A|),K_n(|AA|)\big)$ getting further away from the conjectured limit points. Thus, either $n=32$ is not sufficiently large to see the Sum-Product Conjecture, the Sum-Product Conjecture is false, or we have omitted an unknown important type of set\footnote{For example, our dataset gives an upper bound on the terms of the sequence $\min_{|A|=n} \frac{\log(|A+A|+|AA|)}{\log n}$. Additional work can only shrink these upper bounds. However, the suggestion that Sum-Product Conjecture is false arises from the \emph{trend} in these upper bounds as $n$ grows, not from the upper bounds themselves. More examples could conceivably suggest a different trend.} from our dataset for $n\le 31$.

In Section~\ref{sec:dataset}, we describe the creation of our dataset. In Section~\ref{sec:images}, we show various visualizations of interesting slices of the dataset. In Section~\ref{sec:tiny n}, we prove the exact value of $\SPP(n)$ for $n\le 6$, and describe the obstacle to proving the exact value of $\SPP(7)$. In Section~\ref{sec:questions}, we conclude with a list of problems and conjectures arising from the dataset and its construction. Finally, in Section~\ref{sec:defending} we offer counter-arguments for why the Sum-Product Conjecture remains plausible even without favorable evidence in our dataset.

\section{Constructing the Dataset}
\label{sec:dataset}
The dataset is available at \url{http://www.math.csi.cuny.edu/obryant/SumProduct}. This dataset contains for each $i,j,n$, a set of $n$ positive integers with the smallest known maximum and with sumset having exactly $i$ elements, and product set having exactly $j$ elements. The data is limited to $n\le 32$, and is certainly incomplete even for those $n$.

\begin{figure}
  \includegraphics[width=\textwidth]{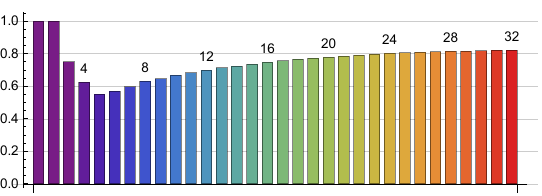}
  \caption{The ratio of the number of pairs $(i,j)$ in $\SPP(n)$ to the number of pairs $(i,j)$ in $[2n-1,n(n+1)/2]^2$. The bars corresponding to $n$ being a multiple of $4$ are labeled. The coloring is consistent with that used throughout this work: the rainbow covers from $n=3$ (purple) to $n=32$ (red).}
  \label{fig:coverage}
\end{figure}

As noted above $\abs{\SPP(n)} \le (n^2-3n+4)^2/4$ for obvious reasons. The fraction of that number that is represented in the dataset is shown in Figure~\ref{fig:coverage}, showing that the dataset contains at least $80\%$ of the truth. In light of the improved bound on $\abs{\SPP(n)}$ given in Corollary~\ref{cor:|SPP(n)|}, we have at least $87\%$ of $\SPP(n)$ for $20\le n \le 32$, and plausibly close to $100\%$.

For $n\le 6$ we prove that the dataset is complete in Section~\ref{sec:tiny n}. There are currently a total of \totalsize\ sum-product pairs in the dataset.

The dataset was obtained by examining
\begin{itemize}[noitemsep]
  \item every subset of $[36] \coloneqq \{1,\dots,36\}$;
  \item every translation and dilation of every set of positive integers with diameter at most $31$;
  \item every subset of the set of divisors of $N$ for all $N\le 2^{11}$ except $1260$, $1440$, $1680$, $1800$, $1980$, and $2016$;
  \item hundreds of pseudorandom subsets of $[N]$ with $n$ elements for each $n \le 32$ and each $N \le 2^{11}$;
  \item hundreds of random subsets of the set of divisors of $N$ for each $N \le 2^{16}$;
  \item the $n$ smallest divisors of each divisor of $2^{11} 3^6 5^4 7^3 11^3 13^2$, for $n\le 32$;
  \item shifts of the sets in the dataset by up to $256$;
  \item For each $A$ in the dataset, we also examined $A\cup\{b\}$ for all $b\in[1,256]$.
\end{itemize}
There were numerous other \emph{ad hoc} attempts that generated few new points, and aren't worth mentioning even though their meager contributions to the dataset remain in the dataset. For example, subsets of geometric progressions and unions of geometric progressions were sampled, but not methodically.

The bulk of the exploration happened with \emph{Mathematica 14}. When speed was critical, {\tt Julia 1.10} was used to isolate the sets that were potentially new, and those were then processed with \emph{Mathematica}. All images were built with \emph{Mathematica}.

\section{Visualizing Sum-Product Pairs}
\label{sec:images}

\subsection{Tiny \texorpdfstring{$n$}{n}}
\begin{figure}
  \begin{center}
  \includegraphics[width=0.85\textwidth]{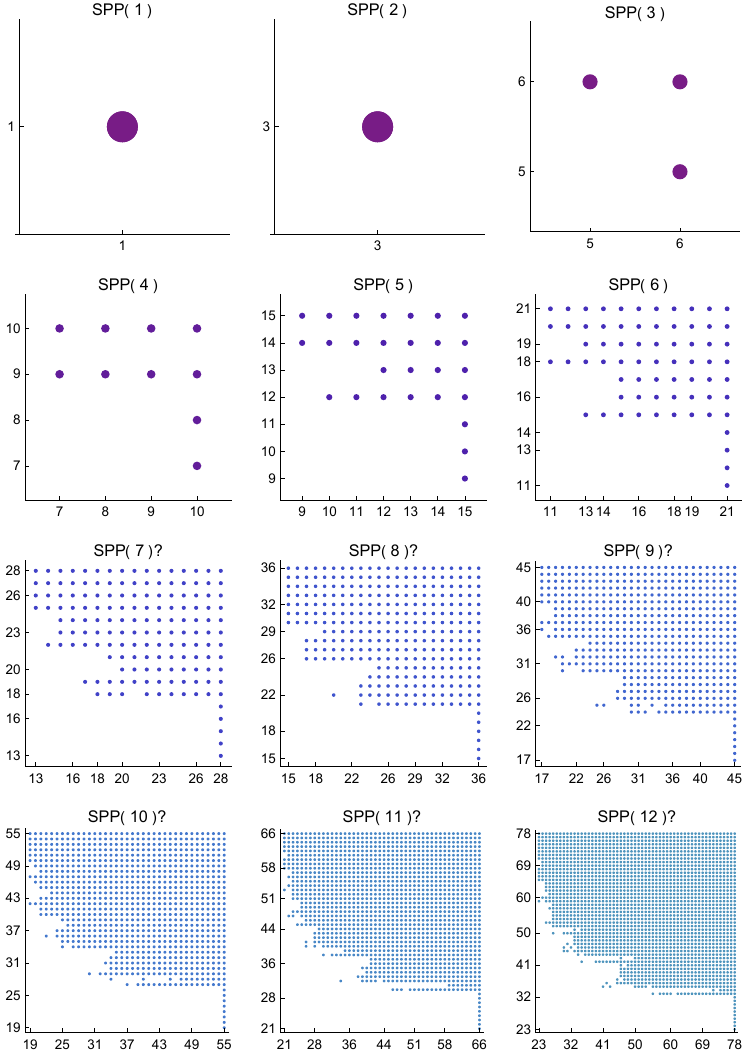}
  \end{center}
  \caption{Plots of $\SPP(n)$ for $n\le 6$, and a large subset of $\SPP(n)$ for $7\le n \le 12$.}
  \label{fig:tiny n}
\end{figure}
We begin with pictures showing the largest subset of $\SPP(n)$ that we have found for $n\le 12$. See Figure~\ref{fig:tiny n}. Note that in the titles of the $12$ pictures, $\SPP(n)$ is followed by a question mark for $n\ge 7$. In Section~\ref{sec:tiny n} we prove that the shown images for $n\le 6$ are the complete $\SPP(n)$.

In each graph, the size of the sumsets is shown on the horizontal abscissae, and the size of the product sets is shown on the vertical ordinates. Arithmetic progressions produce dots along the western edge (on the left); geometric progressions produce dots along the southern edge; Sidon sets\footnote{The set $A=\{a_1<a_2<\dots<a_n\}$ is an arithmetic progression of $a_i=a_1+(i-1)(a_2-a_1)$, a geometric progression if $a_i=a_1(a_2/a_1)^{i-1}$, a Sidon set if $|A+A|=n(n+1)/2$, and a multiplicative Sidon set if $|AA|=n(n+1)/2$.} produce dots along the eastern edge; and multiplicative Sidon sets produce dots along the northern edge.

It is impossible to see the Sum-Product Conjecture in these images\,---\,first because $n$ is so small, but even for large $n$ because almost all integers between $2n-1$ and $n(n+1)/2$ have size roughly $cn^2$. Our normalized pictures below are more informative about the Sum-Product Conjecture.

The first isolated point is $(20,22)\in\SPP(8)$. This point is indeed special. There is only one set of $8$ relatively prime positive integers with exactly $20$ sums and exactly $22$ products: $\{1,2,3,4,6,8,9,12\}$. We do not know if the eight surrounding coordinates are actually not in $\SPP(8)$, but we suspect that this is the case.

We note that $\Psi_8^3 :=\{1,2,3,4,6,8,9,12\}$ is the set of the first 8 numbers all of whose prime factors are at most $3$. Such sets give many extremal examples, so we define some notation for them. We set
  \(\Psi_n^y\)
to the be the set of the first $n$ natural numbers that are $y$-friable (also called $y$-smooth).

\subsection{The {S}idon {E}xclusion {Z}one (SEZ)}
\noindent\begin{minipage}[t]{0.58\textwidth}
\hspace{1em} One quickly notices in Figure~\ref{fig:tiny n}~a flagpole structure. The pole runs down the right edge of each plot, generated by Sidon sets. The rectangular flag waves to the left, tattered along the far edge. This author has no understanding of the tatters, but the flagpole is an actual phenomenon. Specifically, if the product set $AA$ is too small, then $A$ must be a Sidon set.

\hspace{1em} We call this the Sidon Exclusion Zone, abbreviated as ``SEZ'' in our images and throughout this article. This is proved in~\cite{Rice}, and relies on Freiman's $(3n-4)$-Theorem~\cite{Freiman,Freiman2} in the reals.
\end{minipage}
\begin{minipage}[t]{0.40\textwidth}
  \begin{center}
  \vspace{-3mm}
  \includegraphics[width=0.7\textwidth]{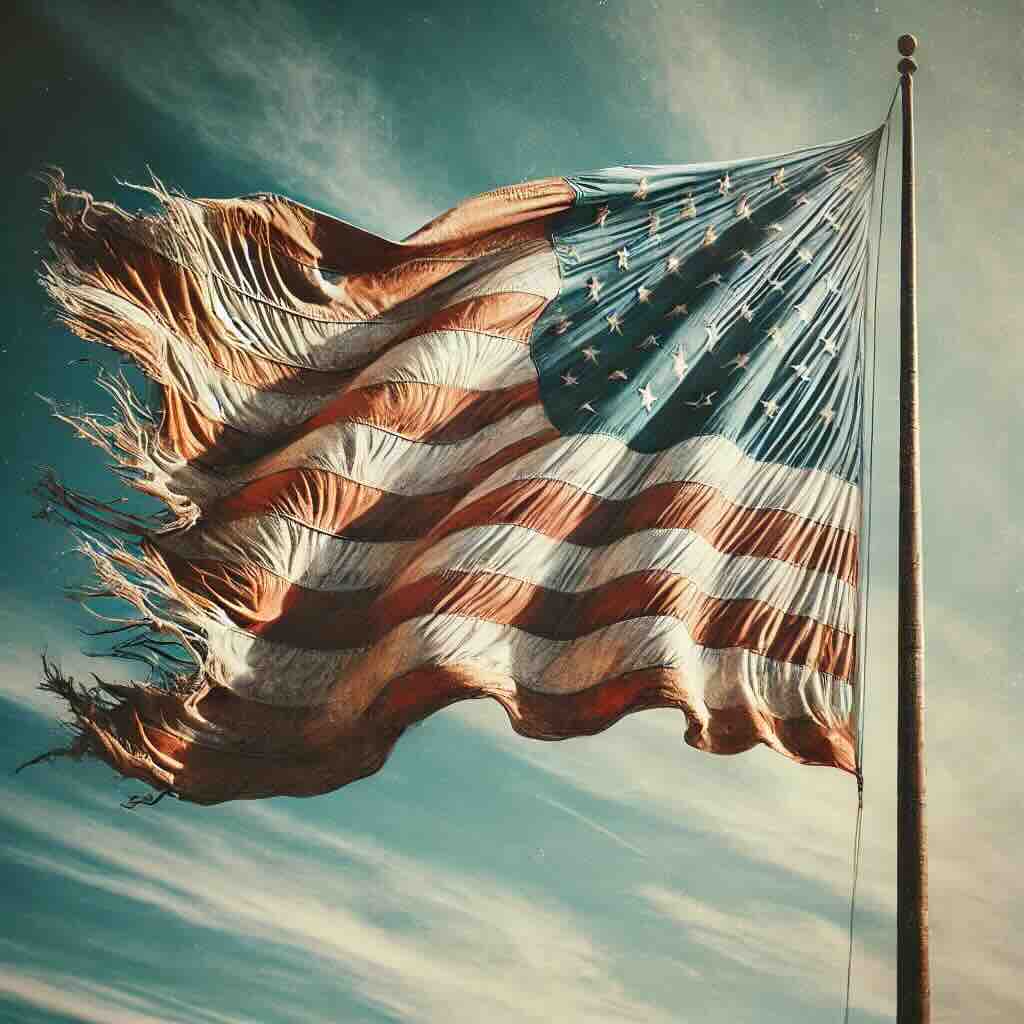}
  \end{center}
  \hspace{3em}\emph{Image created by DALL-E.}
\end{minipage}

\begin{theorem}[Sidon Exclusion Zone~\cite{Rice}]\label{thm:SEZ}
  If $A$ is a set of $n$ positive integers and $|AA| \le 3n-4$, then $A$ is a Sidon set.
\end{theorem}

\begin{theorem}[Freiman's $(3n-4)$-Theorem]
  If $A$ is a set of $n$ real numbers and $|A+A| \le 3n-4$, then $A$ is a subset of an arithmetic progression of length $|A+A|-n+1$.

  If $A$ is a set of $n$ positive real numbers and $|AA| \le 3n-4$, then $A$ is a subset of a geometric progression of length $|AA|-n+1$.
\end{theorem}

A remarkable feature of Theorem~\ref{thm:SEZ} is that we cannot replace ``integers'' with ``reals''. This can be seen by considering geometric progressions like $\{1,\varphi,\varphi^2,\varphi^3\}$ and $\{1,\alpha,\alpha^2,\alpha^3,\alpha^4\}$, where $\varphi^2=\varphi+1$ and $\alpha^3=\alpha+1$.

We can use the Sidon Exclusion Zone to improve the ``easy'' upper bound on $\abs{\SPP(n)}$.

\begin{corollary}\label{cor:|SPP(n)|}
  For $n\ge 3$, we have $\displaystyle \abs{\SPP(n)} \le \frac{(n^2-3n+4)^2}4-\frac{(n-2)^2 (n-1)}{2} $.
\end{corollary}

This author believes that $\abs{\SPP(n)}=\frac 14 n^4 +O(n^3)$ is true, but that $\abs{\SPP(n)}=\frac 14 n^4 -2n^3+O(n^2)$ is false. Perhaps $\abs{\SPP(n)}=\frac 14 n^4 -4n^3+O(n^2)$.

\subsection{Normalizing Correctly}
For $|A|=n$, we know that $|A+A|$ is trapped between $2n-1$ (for arithmetic progressions) and $n(n+1)/2$ (for Sidon sets), inclusive. Also, $2n-1 \le |AA| \le n(n+1)/2$, with the extremal sets being named geometric progressions and multiplicative Sidon sets. The Sum-Product Conjecture draws attention to the exponent of $n$, which ranges from $1$ to $2$ as $x$ ranges from $2n-1$ to $n(n+1)/2$. For small $n$, however, we would like a better normalization than the simple base-$n$ logarithm $L_n(x)=\log_n x$.

We let $\phi_n(x)=m_n x+b_n$. We can choose $m_n,b_n$ in the scenario
\[
  K_n(x) = L_n(x)+\phi_n(x)
\]
so that $K_n(2n-1)=1$ and $K_n(n(n+1)/2)=2$.
Namely,
\begin{align*}
  m_n &= \binom{n-1}{2}^{-1}\,\log _n\left(\frac{4 n-2}{n+1}\right)\\
  b_n &= \frac{n^2-7 n+4}{n^2-3 n+2}\\
  &\qquad +\binom{n-1}{2}^{-1} \left((2n-1) \log _n\left(\binom{n+1}{2}\right)-\binom{n+1}{2} \log_n(2n-1)\right).
\end{align*}
By easy calculus, as $x\le n(n+1)/2$, we have $b_n\to0$ and $m_n x \to 0$.
Alternatively, we can choose $m_n,b_n$ in the scenario
\[ K_n^{(2)}(x) = L_n( \phi_n(x)) = \log_n\left(\frac{2 x - 3n}{1-2/n}\right)\]
so that $K_n^{(2)}(2n-1)=1$ and $K_n^{(2)}(n(n+1)/2)=2$.
A third alternative is to work with
\[ K_n^{(3)}(x) = \phi_n(L_n(x)) = \frac{\log(n(n+1)x)-\log(2(1-2n)^2)}{\log(n(n+1))-\log(4n-2)}=\frac{\log \left(\frac{n (n+1) x}{2 (1-2 n)^2}\right)}{\log \left(\frac{n (n+1)}{4 n-2}\right)},\]
again choosing $\phi_n$ so that  $K_n^{(3)}(2n-1)=1$ and $K_n^{(3)}(n(n+1)/2)=2$.

We prefer $K_n(x)$ as being easiest to use, while $L_n(x)$ is more closely connected to the underlying combinatorics. We include the others primarily to indicate that the phenomena we visualize are in the data and not in the normalization.

\subsection{The Void of Solymosi, and a Result of Bourgain and Chu}
In 2009, Solymosi~\cite{Solymosi} proved that
  \[|A+A|^2 \cdot |AA| \ge \frac{n^4}{1+\floor{\log_2 n}}\]
for any set $A$ of $n$ positive real numbers\footnote{Solymosi states the bound as $n^4/\ceiling{\log n}$, but he has a pair of off-by-one errors in his proof that cancel perfectly, except when $n$ is a power of $2$.}. His elegant combinatorial proof uses two inequalities connecting additive energy to $|A+A|$ and to $|AA|$.

Solymosi's inequality gives
  \[ \log_n \max\{|A+A|,|AA| \} \ge \frac 43 - O\left(\frac{\log\log n}{\log n}\right),\]
and for this reason is often called Solymosi's $4/3$-bound. As $n\to\infty$ it is nearly the strongest result known for $\log_n \max\{|A+A|,|AA| \}$. It has been tweaked more than a half-dozen times since 2009, and the current record is given by Bloom~\cite{Bloom}, improving Solymosi's $4/3$ to $\frac 43+ \frac{2}{951}$ (and allowing $A$ to contain nonpositive numbers).

In terms of our normalization, Solymosi's inequality forbids the limit points of $\NSPP([3,\infty))$ from lying in the region $2x+y \le 4$. We show this region in our pictures with the label ``{Solymosi Void}''.

If $A_n$ is a sequence of multidimensional arithmetic progressions with bounded dimension and cardinality $n$, then $|A+A|\le Cn$ for some constant $C$ depending only on the dimension, and so $K_n(|A+A|)\to 1$. But by the Solymosi Void, the only possible limit point with $x=1$ is $x=K_n(|AA|)=2$. Ergo, if a sequence of sets has linear sumset, it must have quadratic product set (ignoring logarithmic factors).

In light of Theorem~\ref{thm:SEZ} below, we can restrict our attention to sets with $|AA|\ge 3n-3$. Solymosi's inequality is weaker than $|AA| \ge 3n-3,|A+A|\ge 2n-1$ for $n\le 87$, and so it isn't directly relevant to most of our images\footnote{Figure~\ref{fig:smooth numbers} is the exception.}, which are constrained to $n\le 32$. Nevertheless, we include it for visual reference.

The other extreme is to consider $|AA| < C n$, where $A$ is a set of $n$ positive integers, and Chang~\cite{2003.Chang} proved in this case that
\[|A+A| \ge 36^{-C} n^2.\] This result is weaker than $|AA| \ge 3n-3,|A+A|\ge 2n-1$ for $n\le 93\,311$.

\subsection{A Full Description of Figure~\ref{fig:everythingK}}
\label{subsec:everything pics}
Figure~\ref{fig:everythingK} portrays $\NSPP([3,32])$. Figure~\ref{fig:everythingzoom} is a zoomed in version of the image in Figure~\ref{fig:everythingK}. We have first put down a point in red at each normalized sum-product pair from a set with $n=32$ elements. Then, a less red point at each pair from a set with $n=31$, and so on through the rainbow, putting violet points for each pair from a set with $n=3$. Because the points are put down in this order, one doesn't see the mass of points from $n=32$, only the fringe along the bottom edge of the blob of points that isn't covered by $n=31,30,29,\dots$. If we put the points down in the opposite order, one only sees the points arising from $n=32$.

To repeat, red means $n=32$, violet means $n=3$, and the colors are interpolated by the rainbow.

Orienting ourselves in Figure~\ref{fig:everythingK}, the southern edge is for geometric progressions, the western edge for arithmetic progressions, the eastern edge for Sidon sets, and the northern edge for multiplicative Sidon sets.

The northern edge is populated by uniformly random subsets of $N$ for large $N$. That is, the sets are sparse but uniform. The eastern edge comes from random subsets of the set of divisors of $N$ for $N$ up to $2^{16}$. Considering all subsets of the set of divisors for $N\le 2^{11}$ gives the bulk of the lower middle points. Sets along the lower envelope are almost all ``divisor closed'', i.e., if the set contains $x$ then it contains every divisor of $x$. They also tend to almost be $\Psi_n^y$ (the set of the first $n$ positive integers that are $y$-smooth), but frequently with a little deviation.

Figure~\ref{fig:everythingK} (and the other normalizations shown in Figure~\ref{fig:everything4}) casts quantitative doubt on the Sum-Product Conjecture. At the very least, it suggests that $n_0$ will need to be quite large even for modest $\epsilon$, such as $\epsilon=1/2$. However, while the current best results rely on $n$ being quite large, there is a significant chunk of our picture that is empty for small $n$, too.

There are only $28$ sets with $K_n(|A+A|)+2K_n(|AA|) \le 4$, and we enumerate them in Table~\ref{tab:envelope}.
\begin{table}
\[\begin{array}{ccccl}
  \lcm(A) & n & |A+A| & |AA| & \qquad \text{witness} \\ \hline
 288 & 12 & 41 & 35 & \Psi _{11}^3\cup \{32\} \\
144 & 12 & 43 & 35 & \Psi _{11}^3\cup \{36\} \\
864 & 14 & 55 & 43 & \Psi _{14}^3 \\
288 & 14 & 56 & 42 & \Psi _{11}^3\cup \{32,36,48\} \\
864 & 15 & 62 & 46 & \Psi _{15}^3 \\
576 & 15 & 66 & 45 & \Psi _{11}^3\cup \{32,36,48,64\} \\
864 & 16 & 72 & 50 & \Psi _{16}^3 \\
1728 & 16 & 73 & 49 & \Psi _{15}^3\cup \{64\} \\
864 & 17 & 81 & 53 & \Psi _{16}^3\cup \{72\} \\
1728 & 18 & 90 & 57 & \Psi _{18}^3 \\
1728 & 19 & 100 & 60 & \Psi _{18}^3\cup \{96\} \\
3456 & 20 & 115 & 63 & \Psi _{18}^3\cup \{96,128\} \\
21600 & 26 & 101 & 121 & \Psi _{26}^5 \\
21600 & 27 & 108 & 125 & \Psi _{26}^5\cup \{72\} \\
43200 & 28 & 115 & 132 & \Psi _{28}^5 \\
21600 & 28 & 118 & 131 & \Psi _{26}^5\cup \{72,90\} \\
43200 & 29 & 124 & 137 & \Psi _{28}^5\cup \{80\} \\
21600 & 29 & 126 & 136 & \Psi _{26}^5\cup \{72,75,90\} \\
43200 & 30 & 133 & 144 & \Psi _{28}^5\cup \{80,90\} \\
43200 & 30 & 136 & 141 & \Psi _{28}^5\cup \{80,96\} \\
43200 & 31 & 140 & 150 & \Psi _{30}^5\cup \{90\} \\
43200 & 31 & 142 & 149 & \Psi _{28}^5\cup \{80,90,100\} \\
43200 & 32 & 150 & 155 & \Psi _{30}^5\cup \{90,96\} \\
129600 & 33 & 156 & 163 & \Psi _{33}^5 \\
129600 & 35 & 175 & 174 & \Psi _{35}^5 \\
129600 & 36 & 183 & 178 & \Psi _{36}^5 \\
1296000 & 47 & 313 & 245 & \Psi _{47}^5 \\
\end{array}\]
  \caption{The 28 known sets whose normalized sum-product pair is in the region $x+2y\le 4$.}
  \label{tab:envelope}
\end{table}

Conjecture~\ref{conj:SV} is motivated by the paucity of sets that are even close to being counterexamples.

\begin{figure}
  \begin{center}
    \includegraphics[width=\textwidth]{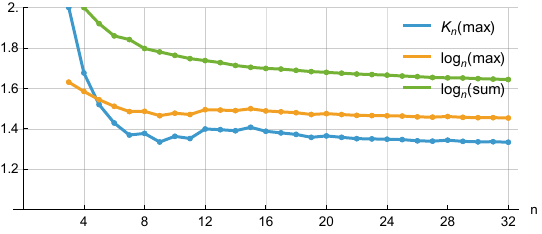}
  \end{center}
  \caption{The minimum values of $K_n(\max\{|A+A|,|AA| \})$, $\log_n(\max\{|A+A|,|AA| \})$, and $\log_n(|A+A|+|AA|)$ in the dataset, as a function of $n$. The three functions converge to each other, and the Sum-Product Conjecture says that they all converge to 2.}
  \label{fig:SPnumbers}
\end{figure}

Another interesting chart is the minimum value of
  \[K_n(\max\{|A+A|,|AA| \})\]
that we have found, as a function of $|A|$. Or one may prefer to consider
  \[\min_A \log_n(\max \{|A+A|,|AA| \}).\]
Or perhaps one feels the proper quantity to consider is
  \[\min_A \log_n(|A+A|+|AA|).\]
These are all shown in Figure~\ref{fig:SPnumbers}. The Sum-Product Conjecture implies that all 3 of these quantities should go to 2 as $n\to\infty$. In the image, not only do we not see a trend to 2, they are steadily decreasing.

Based on Figure~\ref{fig:SPnumbers}, we make Conjecture~\ref{conj:|A+A|+|AA|}, stated in the Introduction.
We confess that the appearance of the golden ratio in Conjecture~\ref{conj:|A+A|+|AA|} is whimsical, but nevertheless is plausible given the data. It cannot be replaced with $1.645$, as witnessed by $\Psi_{29}^7\cup \{50,54,60\}$.
This set has $32$ elements, $104$ sums, and $194$ products.

In Table~\ref{tab:optimals} we give for each $n$ a set of $n$ positive integers that minimizes $\max\{|A+A|,|AA|\}$ in our dataset.
\begin{table}
  \[\begin{array}{cccc}
    n & |A+A| & |AA| & \text{optimal set} \\ \hline
 1 & 1 & 1 & \Psi _1^2 \\
2 & 3 & 3 & \Psi _2^2 \\
3 & 6 & 5 & \Psi _3^2 \\
4 & 7 & 9 & \Psi _4^3 \\
5 & 10 & 12 & \Psi _5^3 \\
6 & 13 & 15 & \Psi _6^3 \\
7 & 18 & 18 & \Psi _6^3\cup \{12\} \\
8 & 20 & 22 & \Psi _8^3 \\
9 & 25 & 25 & \Psi _9^3 \\
10 & 30 & 29 & \Psi _{10}^3 \\
11 & 34 & 32 & \Psi _{11}^3 \\
12 & 41 & 35 & \Psi _{11}^3 \cup \{32\} \\
13 & 43 & 46 & \Psi _7^5 \cup \{10,12,16,20,24,32\} \\
14 & 51 & 50 & \Psi _7^5 \cup \{10,12,16,20,24,32,40\} \\
15 & 45 & 58 & \Psi _{11}^3 \cup \{5,10,20,32\}\\
16 & 60 & 62 & \Psi _{11}^3 \cup  \{10,20,32,36,48 \} \\
17 & 61 & 67 & \Psi _{11}^3 \cup  \{5,10,20,32,36,48 \} \\
18 & 64 & 72 & \Psi _{23}^5 \setminus  \{15,25,27,30,45 \}\\
19 & 75 & 76 & \Psi _{11}^3 \cup  \{5,10,20,32,36,40,48,64 \} \\
20 & 81 & 83 & \Psi _{15}^3 \cup  \{5,10,15,20,64 \}\\
  21 & 83 & 88 & \Psi_{15}^5 \cup  \{30,32,36,40,48,64 \} \\
22 &91 &93 & \Psi_{28}^5 \setminus  \{25,27,45,50,54,64 \} \\
23 &97 &99 & \Psi_{28}^5 \setminus  \{25,27,45,50,54 \} \\
24 &98 &105 & \Psi_{28}^5 \setminus  \{25,45,50,64 \} \\
25 &102 &111 & \Psi_{28}^5 \setminus  \{25,50,64 \} \\
26 &114 &116 & \Psi_{32}^5 \setminus  \{25,50,64,75,80,81 \} \\
27 &122 &121 & \Psi_{33}^5 \setminus  \{25,50,75,80,81,90 \} \\
28 &124 &130 & \Psi_{23}^5 \cup  \{54,60,64,72,96 \} \\
29 &134 &135 & \Psi_{26}^5 \cup  \{72,90,108 \} \\
30 &136 &141 & \Psi_{28}^5 \cup  \{80,96 \} \\
31 &148 &147 & \Psi_{35}^5 \setminus  \{64,80,96,100 \} \\
32 &154 &154 & \Psi_{35}^5 \setminus  \{75,81,100 \}
  \end{array}
  \]
  \caption{For each $n$, a set that minimizes $\max\{|A+A|,|AA|\}$ among sets with $n$ elements. These are not known to be optimal sum-product pairs (and in many cases are not the unique sum-product pair fitting the description).}
  \label{tab:optimals}
\end{table}

In another we direction, we conjecture a bound on the maximum size of a smallest example in $\SPP(n)$.
\begin{conjecture}
  If $(i,j)\in \SPP(n)$, then there is a set $A\subseteq\NN$ with $|A|=n$, $|A+A|=i$, $|AA|=j$, and $\max A \le 2^{\floor{(3n-4)/2}}$.
\end{conjecture}
This maximum seems to happen at the top of the flagpole ($|AA|=3n-4,|A+A|=n(n+1)/2$), but we are unable to prove that there isn't some sum-product pair elsewhere that requires an even higher maximum.

\subsection{Special Categories of Sets}
\begin{figure}
  \begin{center}
    \includegraphics[width=\textwidth]{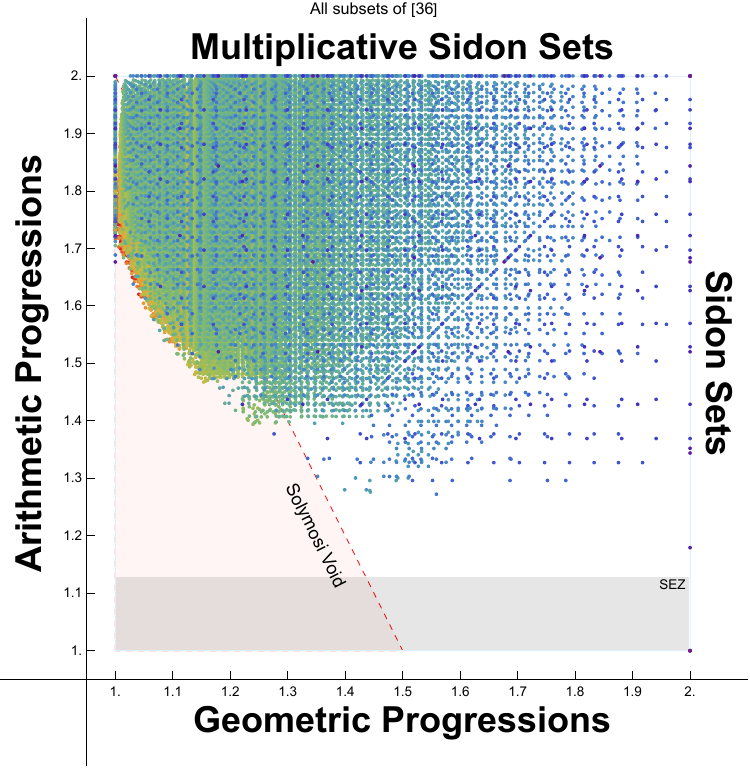}
  \end{center}
  \caption{The points $\big(K_n(|A+A|),K_n(|AA|)\big)$ for all $A\subseteq[36]$.  Redder points come from larger $|A|$.
    }
  \label{fig:AP36}
\end{figure}

The most obvious way to start populating a dataset with small sets with surprising values of $(|A+A|,|AA|)$ to loop through all subsets of $[N]\coloneqq\{1,2,\dots,N\}$. We have done this for $N=36$ on our desktop. We estimate that an excellent programmer with access to a nice cluster could examine every subset of $[45]$, or perhaps even $[50]$. Figure~\ref{fig:AP36} shows the result of our computation.

There are several visual artifacts that are not meaningful. One sees lines sloping from southwest to northeast, and more graceful arcs bowed away from the origin from the north edge to east edge. These are both the result of the grid of possible values being normalized. There are light vertical streaks around $1.22$ and $1.12$; these are the result of the irregular distribution resulting from the $32$ normalized lattices being combined.

One sees that as $n$ grows, the dots are not pushed north and east, as suggested by the Sum-Product Conjecture, but rather are pushed mostly west and slightly south. The western push is a red-herring, as $A\subseteq[36]$ implies that $A+A\subseteq[2,72]$, which strongly restricts $A$ from having a large sumset relative to diameter for $|A|\ge 20$. We believe the slight southern drift is an effect of small numbers, in particular the large number of solutions to $ab=cd$ among $4$-tuples from $[36]$, which are difficult to avoid when taking a large subset of $[36]$.

\begin{figure}
  \begin{center}
    \includegraphics[width=0.49\textwidth]{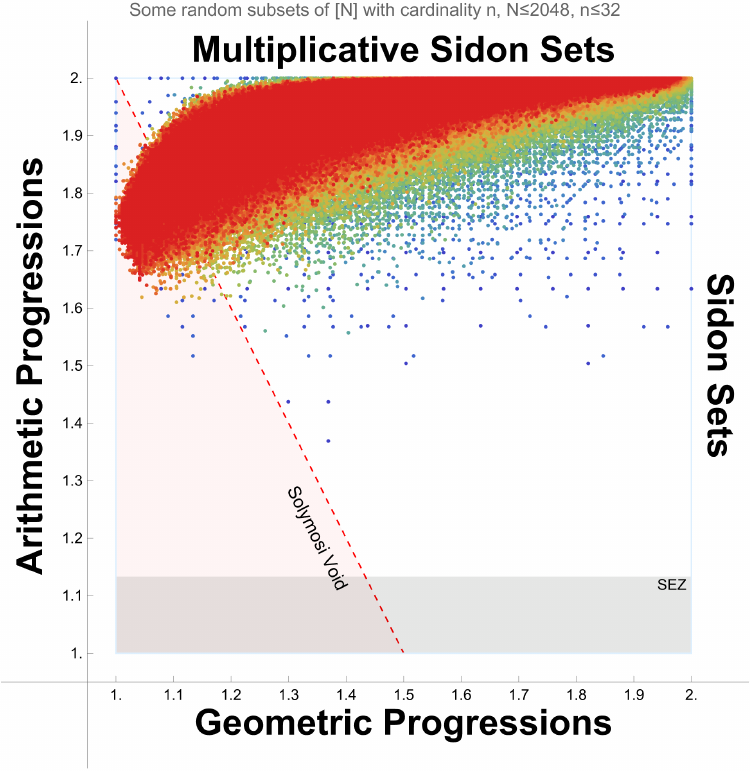}
  \end{center}
  \caption{The points $\big(K_n(|A+A|),K_n(|AA|)\big)$ for a large number of uniformly random subsets of $N$ with cardinality $n$, where $33\le N \le 2048$ and $5\le n \le 32$.}
  \label{fig:randomsubsets}
\end{figure}

We have generated hundreds of pseudorandom subsets of $N$ with $n$ elements for $N\le 2^{16}$ and $n\le 32$. When $N$ is small, these tend to be nearly arithmetic progressions (small sumset, large product set), and when $N$ is large they tend to be almost Sidon sets and multiplicative Sidon sets. The image in Figure~\ref{fig:randomsubsets} shows this slice of the dataset. The line $y=x$ is pronounced. Is it more likely for a uniformly random subset of $[N]$ to have $|A+A|=|AA|$ than, say, $|A+A|=|AA|+3$?

\begin{figure}
  \begin{center}
    \includegraphics[width=0.49\textwidth]{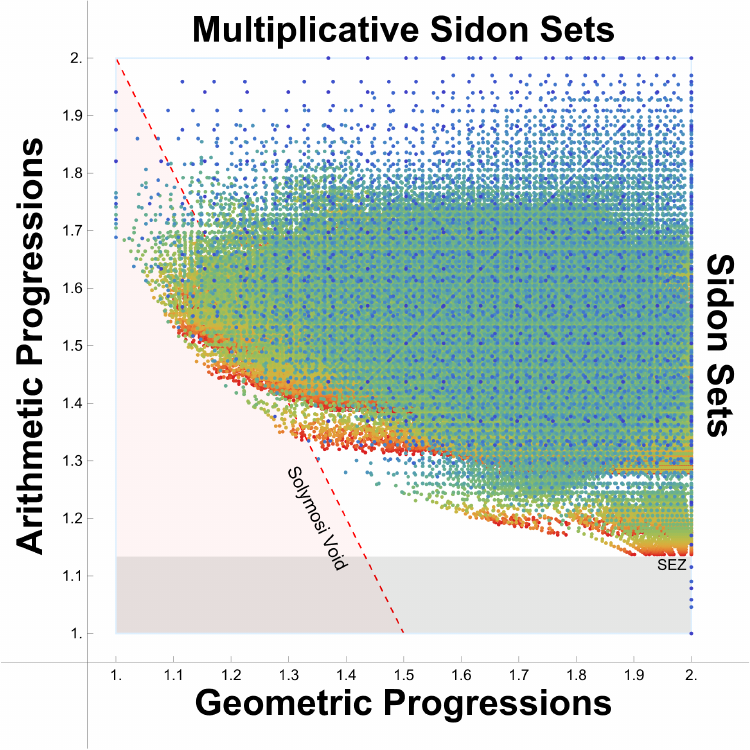}
    \includegraphics[width=0.49\textwidth]{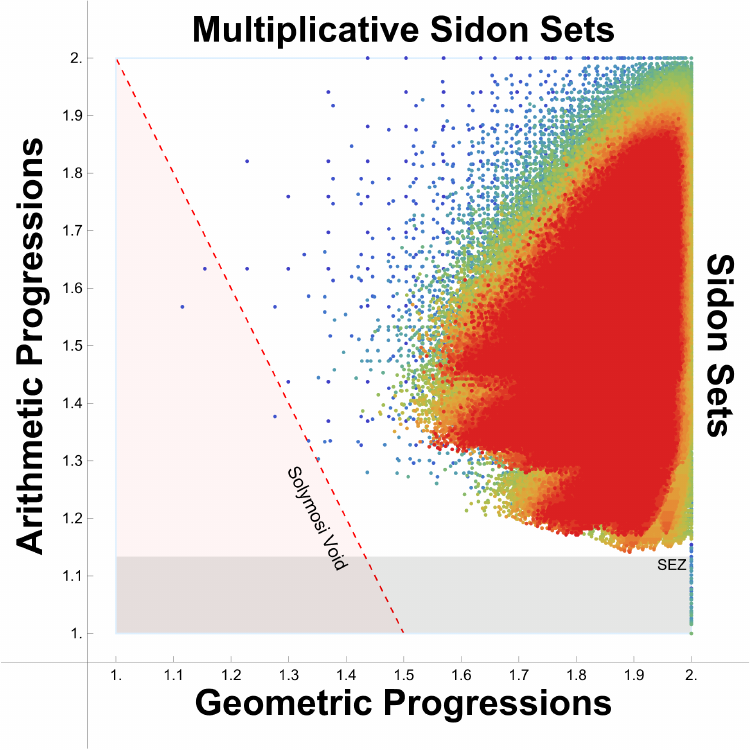}
  \end{center}
  \caption{On the left are the points $\big(K_n(|A+A|),K_n(|AA|)\big)$ for every subset of the set of divisors of $N$, where $N\le 2048$ (except for six values of $N$ mentioned in the text). The image on the right is generated by considering random subsets of the set of divisors of $N$ for all $N\le 65536$.}
  \label{fig:divisorsets}
\end{figure}

Another interesting category of sets is sets of divisors. These are typically not having sum-product pairs in an interesting location, but the set of proper divisors does much better. The set of divisors of $N$ that are most $\sqrt{N}$ often gets into the middle of our image. Figure~\ref{fig:divisorsets} shows the contribution of all subsets of the divisors of $N$ for $N\le 2048$ except $N\in \{ 1260, 1440, 1680, 1800, 1980,2016\}$.

We also considered a large random sampling of subsets of size $n$ of the set of divisors of $N$ for all $N\le 65536$, with $n\le32$. This image is also shown in Figure~\ref{fig:divisorsets}.

\subsubsection{Sets with Small Sumset or Small Product Set}
\label{subsec:small sumsets}
\begin{theorem}\label{thm:small sumset}
  If $A$ is a set of $n$ positive real numbers with $|A+A|\le 3n-4$ and $|AA|<n(n+1)/2$, then some dilation of $A$ is a subset of
  \[\frac{i j-k \ell}{k+\ell-i-j}+[0,2n-4]\]
  where $0\le i <k \le \ell < j \le 2n-4$. In particular,
  \begin{multline*}
    \SPP_{\RR_+}(n)\cap\left([2n-1,3n-4]\times[2n-1,\tfrac{n(n+1)}2-1]\right)\\
      =\SPP(n)\cap\left([2n-1,3n-4]\times[2n-1,\tfrac{n(n+1)}2-1]\right).
  \end{multline*}
\end{theorem}

\begin{proof}
Suppose that $A$ is a set of $n$ real numbers and $|A+A|\le 3n-4$, so that Freiman's $(3n-4)$-Theorem guarantees that $A$ is contained in an arithmetic progression of length at most $2n-3$. Say $A\subseteq \{a,a+d,\dots,a+(2n-4)d\}$. Suppose further that $|AA|<n(n+1)/2$ strictly, so that some pair of products are equal. Say $(a+id)(a+jd)=(a+kd)(a+\ell d)$ with $0\le i <k \le \ell < j \le 2n-4$. Thus,
  \[a= d\,\cdot \, \frac{i j-k \ell}{k+\ell-i-j}.\]
By dilating if necessary, we can now assume that $A$ is a subset of
  \[\frac{i j-k \ell}{k+\ell-i-j} + [0,2n-4],\]
which has $2n-3$ elements, for some $0\le i <k \le \ell < j \le 2n-4$.
\end{proof}

There are a small number of $\frac{i j-\ell k}{\ell-i-j+k}$, and a small number of subsets of size $n$ of a set of size $2n-3$. At least, for $n\le 10$, this is a tractable search. We have added all of these sets to our dataset. This is noteworthy because it means that our dataset is complete for $n\le 10$ and $|A+A| \le 3n-4$ for sets of reals, and therefore for sets of integers, too.\footnote{I give thanks to Moshe Newman for suggesting this approach.}

For example, with $n=6$, we find that $a$ is among the 23 rationals
\begin{equation}\label{eq:a possibilities}
    a \in \left\{\frac{1}{6},\frac{1}{5},\frac{1}{4},\frac{1}{3},\frac{2}{5},\frac{1}{2},\frac{2}{3},\frac{3}{4},1,\frac{4}{3},\frac{3}{2},2,\frac{5}{2},3,4,\frac{9}{2},5,6,7,8,9,10,12\right\},
\end{equation}
so that any $(i,j)$ in $\SPP_{\RR_+}(6)$ with $i\le 14$ and $j\le20$ is witnessed by $a+[0,2(6)-3]$ for some $a$ in Line~\eqref{eq:a possibilities}. There are just $\binom{10}{6}$ such sets, and $23$ possible values of $a$, giving just $4830$ sets to consider. A large majority of these have a maximum of at most $36$, and so have already been incorporated into our dataset.

\begin{theorem}\label{thm:small product}
  Suppose that $A$ is a set of $n$ positive real numbers with $|AA|\le 3n-4$ and $|A+A|<n(n+1)/2$. Then $A$ is a subset of a geometric progression of length $|A+A|-n+1$ whose ratio is irrational, between 1 and 2, and a root of a polynomial of the form $\sum_{d=\ell}^{j-1} x^d - \sum_{d=0}^{k-1} x^d$, where $0< k \le \ell < j \le |A+A|-n \le 2n-4$.
\end{theorem}

\begin{proof}
If $n\le 2$, the hypotheses are not satisfied, so assume that $n\ge 3$.
By dilation, which affects neither $|A+A|$ nor $|AA|$, we may assume that $\min A=1$. By Frieman's $(3n-4)$-Theorem, we know that $A$ is a subset of a geometric progression with length $|A+A|-n+1$, say
  $A\subseteq G\coloneqq \{1,r,r^2,\dots,r^{|A+A|-n}\}$
where $r>1$.

As $|G|=|A+A|-n\le 2n-4$ and $|A|=n$ and $n\ge 3$, the set $A$ must contain consecutive elements of $G$. If their ratio, $r$, is rational, then $G$ is a set of rational numbers. Consequently, $A$ must be a set of rational numbers, which we can dilate to a set of integers. The SEZ then tells us that $A$ does not exist. Consequently, $r$ is irrational.

Since $|A+A|<n(n+1)/2$ strictly, there must be some $0\le i <k\le \ell < j\le |A+A|-n$ with $r^i+r^j=r^k+r^\ell$. Dividing through by $r^i$, and relabeling $k,j,\ell$ if necessary, we find that $r$ is a root of $p(x)=x^j-x^\ell-x^k+1$ for some $0< k\le \ell < j\le |A+A|-n$. Since $j>\ell\ge k$, we see that $p(x)\ge 1$ for $x\ge 2$, so that $r<2$. Also, $p(1)=0$, whence $r$ is a root of $p(x)/(x-1)= \sum_{d=0}^{k-1} x^d - \sum_{d=\ell}^{j-1} x^d$.
\end{proof}

\subsection{Number Usage}
We believe that our dataset does not contain every sum-product pair even for $n=12$, much less for $n=32$. And among those sum-product pairs that are in the dataset, it is highly unlikely that we found the set with the smallest maximum that has that pair.

\begin{figure}
  \begin{center}
    \includegraphics[width=\textwidth]{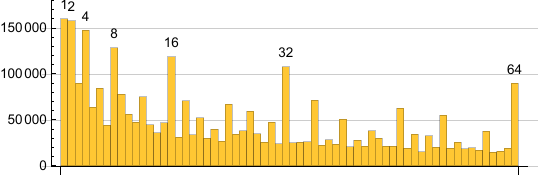}
  \end{center}
  \caption{The number of sets with $n=32$ in our dataset that contain $k$, for $1\le k \le 64$. The bars for the values $k=1,2,4,\dots,64$ are labeled.}
  \label{fig:usage}
\end{figure}

Nevertheless, we were curious as to which numbers contribute more often. Figure~\ref{fig:usage} contains a bar chart of the number of sets in our dataset with $n=32$ that contain the number $k$, $1\le k \le 64$. Even numbers are more likely: $3\,618\,950$ even numbers and just $2\,094\,554$ odd numbers were used with $n=32$.

\subsection{Friable Numbers}
A natural candidate for a set that is multiplicatively small are the friable numbers, also called the smooth numbers. As above, $\Psi_n^y$ is the set of the $n$ smallest positive integers all of whose prime factors are at most $y$.

\begin{figure}
  \begin{center}
  \includegraphics[width=\textwidth]{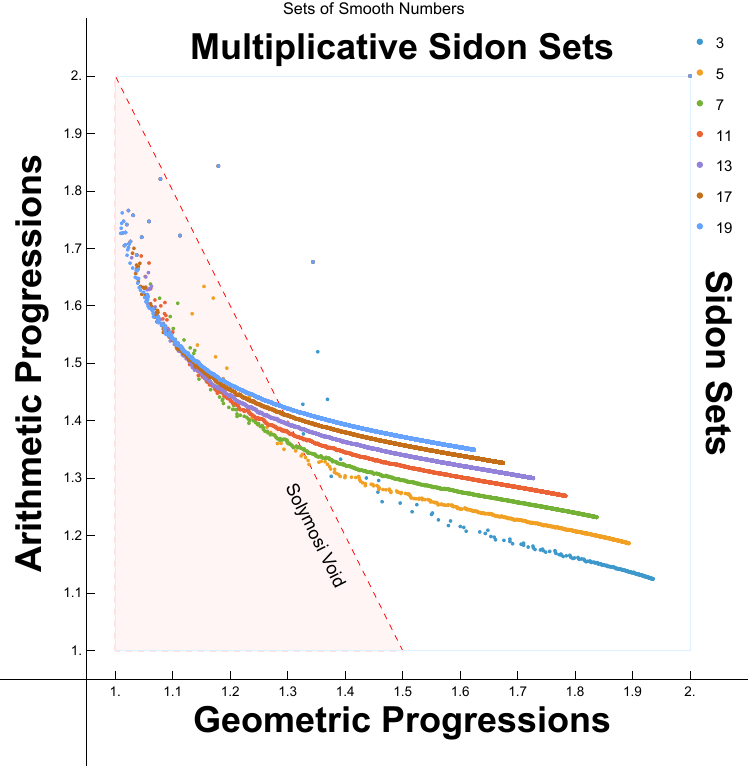}
\end{center}
  \caption{The normalized sum-product pairs for all $\Psi_n^y\subseteq[1,2^{24}]$ for prime $y$ up to $19$. The set $\Psi_n^y$ is the set of the $n$ smallest numbers all of whose prime factors are at most $y$. This image contains $57\,999$ points.}
  \label{fig:smooth numbers}
\end{figure}

In Figure~\ref{fig:smooth numbers} on sees the normalized sum-product pairs for $\Psi_n^y$ for $y\in \{3,5,7,11,13,17,19\}$ and $n\ge3$ with $\max\Psi_n^y\le 2^{24}$. That's $201$ sets that are $3$-friable, and $24\,090$ sets that are $19$-friable. {201, 835, 2401, 5141, 9597, 15748, 24090}

One sees in this image that for fixed $y$, the sumset size $|A+A|$ becomes more an more quadratic as $n$ grows. On the other and, if one fixes a size $|A+A|\approx |A|^x$, then the product set slowly climbs in size towards $|A|^2$. Thus, the sets $\Psi_n^y$ that feature so prominently in many of our best examples are not on their own going to provide a counterexample to the Sum-Product Conjecture.

\section{Sum-Product Pairs for \texorpdfstring{$n \le 6$}{n ≤ 6}}
\label{sec:tiny n}
We define
\[ \SPP_X(n) \coloneqq  \left\{ \big(|A+A|,|AA| \big) : A \subseteq X, |A|=n \right\},\]
where $X$ can be any subset of a ring.
In the most common setting, $X=\NN$, we omit it from our notation.

The purpose of this section is to provide hand proofs of the values of $\SPP(n)$ for $n\le 4$, and to indicate the needed computations to prove the value of $\SPP(n)$ and $\SPP_{\RR_+}(n)$ for $n\le 6$.

\begin{theorem}\label{thm:main}
  The sets $\SPP(n),\SPP_{\RR_+}(n)$ for $n\le 4$ have the following values:
  \begin{align*}
    \SPP_{\RR_+}(1) &= \SPP(1) = \{(1,1)\}\\
    \SPP_{\RR_+}(2) &= \SPP(2) = \{(3,3)\} \\
    \SPP_{\RR_+}(3) & =\SPP(3) = \{(5,6),(6,5),(6,6)\}\\
    \SPP(4) &=\big\{ (i,j) : 7 \le i \le 10, 9 \le j \le 10\big\} \cup\big\{(10,7),(10,8)\big\} \\
    \SPP_{\RR_+}(4) &= \SPP(4) \cup\{(9,7), (9,8)\}
  \end{align*}
\end{theorem}

\begin{proof}
The values of $\SPP(1)$ and $\SPP(2)$ are trivial. The values of $\SPP(3)$ and $\SPP(4)$ follow immediately from the three facts:
$\SPP(n) \subseteq [2n-1,n(n+1)/2]^2$; Theorem~\ref{thm:SEZ}; and the examples:
\begin{center}
  \includegraphics[width=\textwidth]{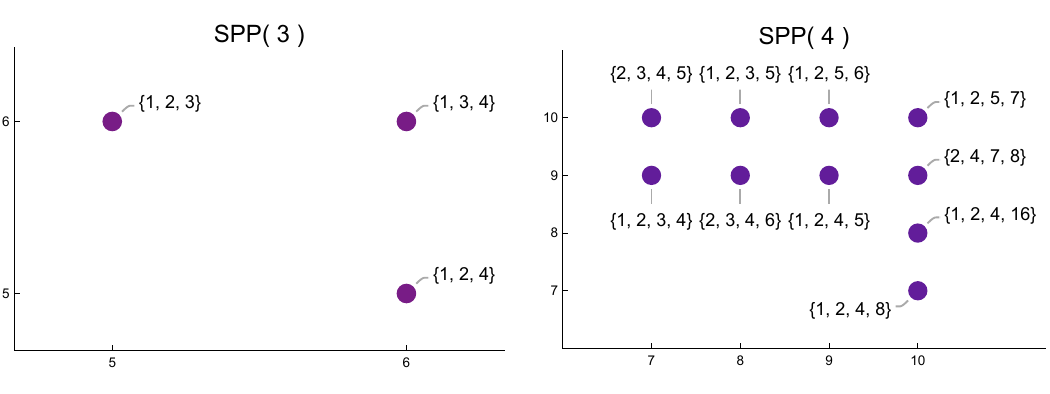}
\end{center}

For $\SPP_{\RR_+}(4)$, we need to identify which of the six potential sum-product pairs $(7,7)$, $(8,7)$, $(9,7)$, $(7,8)$, $(8,8)$, $(9,8)$ can be realized by sets of real numbers. All six potential pairs imply a set $A$ that satisfies the hypotheses of Theorem~\ref{thm:small product}, and so $A$ must be a dilation of a subset of the geometric progression $\{1,r,r^2,r^3,r^4\}$, where $1<r<2$ is irrational and a root of one of $x^2-x-1$, $x^3-x-1$, $x^3-x^2-x-1$. For each of these three real numbers, and for each of the $5$ subsets of with $4$ elements, we compute the sumset and product set. With $r^2-r-1=0$, the set $\{1,r,r^2,r^3\}$ has $9$ sums and $7$ products; the set $\{1,r^2,r^3,r^4\}$ has $9$ sums and $8$ products. No other sum-product pairs arise.
\end{proof}

\begin{theorem}
  The  sets $\SPP(5),\SPP(6)$ are
  \begin{align*}
    \SPP(5) &= \big\{ (i,j) : 12\le i \le 15, 12 \le j \le 15\big\}  \cup \big\{(9,14),(9,15),(10,12),\\
    &\qquad (10,14),(10,15),(11,12),(11,14),(11,15),
     (15,9),(15,10),(15,11) \big\},\\
    \SPP(6) &= \big\{(i,j) : 15\le i \le 21,15\le j \le 21 \big\}
      \cup\{(21,11),(21,12),(21,13),\\
      &\qquad (21,14),(13,15),(14,15) ,(14,19),(13,19) \big\} \\
      &\qquad \cup \{(i,j):11\le i \le 14,j\in\{18,20,21\}\big\}
  \end{align*}
 and the sets $\SPP_{\RR_+}(5),\SPP_{\RR_+}(6)$ are
  \begin{align*}
    \SPP_{\RR_+}(5) &= \SPP(5) \cup\{(i,j): 13\le i \le 14, 9 \le j \le 11\}\\
    \SPP_{\RR_+}(6) &= \SPP(5) \cup\{(i,j): 18\le i \le 20, 11 \le j \le 14\}
  \end{align*}
\end{theorem}

\begin{proof}
First, we give $27$ sets\footnote{As an added bonus, these are sets with $5$ elements and with minimum maximum for each sum-product pair.} that generate each of the $27$ elements of $\SPP(5)$.
\[
\begin{array}{ccc|ccc}
  |A+A| & |AA| & A & |A+A| & |AA| & A \\ \hline
 12 & 12 & \{1,2,3,4,8\} & 9 & 14 & \{1,2,3,4,5\} \\
15 & 12 & \{1,2,4,8,13\} & 12 & 13 & \{1,2,3,4,9\} \\
13 & 12 & \{1,2,4,7,8\} & 12 & 15 & \{1,2,5,6,7\} \\
11 & 15 & \{1,2,3,5,7\} & 15 & 13 & \{1,2,4,11,16\} \\
14 & 12 & \{1,2,4,8,9\} & 15 & 15 & \{1,3,8,11,12\} \\
10 & 14 & \{3,4,5,6,8\} & 13 & 13 & \{1,2,4,6,9\} \\
15 & 11 & \{1,2,4,16,32\} & 13 & 15 & \{1,2,5,6,8\} \\
15 & 9 & \{1,2,4,8,16\} & 11 & 14 & \{1,2,3,5,6\} \\
12 & 14 & \{1,2,3,6,7\} & 14 & 13 & \{1,3,4,6,12\} \\
15 & 14 & \{1,2,5,10,12\} & 14 & 15 & \{1,2,6,7,9\} \\
13 & 14 & \{1,2,4,6,7\} & 15 & 10 & \{1,2,4,8,32\} \\
14 & 14 & \{1,3,4,8,9\} & 10 & 12 & \{1,2,3,4,6\} \\
11 & 12 & \{2,3,4,6,8\} & 9 & 15 & \{3,4,5,6,7\} \\
10 & 15 & \{1,3,4,5,6\} & & &
\end{array}
\]
The SEZ excludes the possibility of $(i,j)$ with $9\le i \le 14$, $9\le j \le 11$. It remains to exclude the four points $(9,12)$, $(9,13)$, $(10,13)$, and $(11,13)$. The search described in Subsection~\ref{subsec:small sumsets} would have found any set giving one of those pairs. This concludes the proof of $\SPP(5)$.

We now consider $\SPP(6)$. In Table~\ref{tab:n=6} we give a witness for each of the $69$ pairs in $\SPP(6)$. The pairs can be seen in Figure~\ref{fig:tiny n}.
\begin{table}
\[
\begin{array}{ccc|ccc}
  |A+A| &|AA| & \text{witness} & |A+A| & |AA| & \text{witness} \\ \hline
  16 & 18 & \{1,2,3,6,7,9\} & 14 & 15 & \{1,2,3,4,6,9\} \\
  13 & 20 & \{1,2,3,5,6,7\} & 17 & 17 & \{1,2,4,6,8,9\} \\
  14 & 20 & \{1,2,3,5,6,8\} & 20 & 15 & \{1,2,4,8,11,16\} \\
  15 & 17 & \{1,2,3,4,8,9\} & 18 & 17 & \{1,2,3,6,8,12\} \\
  20 & 16 & \{1,3,4,6,12,24\} & 21 & 16 & \{3,5,6,12,20,24\} \\
  20 & 20 & \{1,2,5,6,12,14\} & 16 & 15 & \{1,2,3,4,6,12\} \\
  16 & 16 & \{2,3,4,6,9,12\} & 21 & 20 & \{1,2,9,12,14,18\} \\
  16 & 20 & \{1,2,3,5,8,9\} & 21 & 19 & \{1,2,6,9,18,20\} \\
  14 & 19 & \{1,2,3,4,5,9\} & 11 & 18 & \{1,2,3,4,5,6\} \\
  13 & 19 & \{2,4,6,7,8,9\} & 11 & 20 & \{2,3,4,5,6,7\} \\
  19 & 21 & \{2,3,5,9,11,12\} & 17 & 15 & \{1,2,4,5,8,10\} \\
  20 & 19 & \{1,2,4,5,11,16\} & 18 & 15 & \{1,2,3,4,8,16\} \\
  16 & 19 & \{1,2,4,6,7,9\} & 21 & 21 & \{1,2,5,11,13,18\} \\
  17 & 18 & \{1,2,3,6,9,10\} & 12 & 18 & \{2,3,4,5,6,8\} \\
  18 & 18 & \{1,2,3,6,9,11\} & 15 & 18 & \{1,2,4,5,7,8\} \\
  15 & 15 & \{2,3,4,6,8,12\} & 12 & 20 & \{1,2,3,4,5,7\} \\
  17 & 16 & \{1,2,3,4,8,12\} & 15 & 16 & \{2,3,4,6,8,9\} \\
  17 & 20 & \{1,2,3,7,9,10\} & 15 & 20 & \{1,2,3,6,7,8\} \\
  18 & 16 & \{2,3,4,6,9,18\} & 21 & 17 & \{1,2,4,11,16,22\} \\
  18 & 20 & \{1,2,4,5,9,11\} & 11 & 21 & \{5,6,7,8,9,10\} \\
  19 & 17 & \{1,2,4,5,10,16\} & 15 & 19 & \{1,2,3,4,7,9\} \\
  14 & 21 & \{1,2,5,6,7,8\} & 19 & 18 & \{1,2,4,8,10,11\} \\
  13 & 21 & \{1,3,5,6,7,8\} & 19 & 16 & \{1,2,3,6,12,18\} \\
  20 & 21 & \{1,2,6,9,13,15\} & 19 & 20 & \{1,2,5,7,10,11\} \\
  17 & 19 & \{1,2,4,6,9,10\} & 21 & 15 & \{1,2,4,8,16,21\} \\
  16 & 21 & \{1,2,5,6,8,9\} & 12 & 21 & \{1,3,4,5,6,7\} \\
  18 & 19 & \{1,2,3,8,9,12\} & 15 & 21 & \{1,2,3,5,7,8\} \\
  21 & 11 & \{1,2,4,8,16,32\} & 19 & 19 & \{2,3,6,8,11,12\} \\
  20 & 17 & \{1,3,4,6,12,16\} & 13 & 18 & \{1,2,3,4,6,7\} \\
  19 & 15 & \{1,2,4,7,8,16\} & 14 & 18 & \{1,2,3,4,7,8\} \\
  16 & 17 & \{2,3,6,8,9,12\} & 20 & 18 & \{1,3,4,11,12,16\} \\
  17 & 21 & \{1,2,3,7,8,10\} & 21 & 12 & \{1,2,4,8,16,64\} \\
  18 & 21 & \{1,2,6,8,9,11\} & 21 & 13 & \{1,2,4,8,32,64\} \\
  21 & 18 & \{1,2,4,8,13,21\} & 21 & 14 & \{1,2,4,8,64,128\} \\
  13 & 15 & \{1,2,3,4,6,8\} & & &
\end{array}
\]
\caption{A set with $6$ elements giving each sum-product pair possible. Those witnesses with maximum at most 36 have the smallest maximum possible.}
\label{tab:n=6}
\end{table}
The SEZ excludes the possibility of $(i,j)$ with $11\le i \le 20$, $11\le j \le 14$. The search described in Subsection~\ref{subsec:small sumsets} found all pairs with $11\le i \le 14,11\le j \le 20$. This concludes the proof of $\SPP(6)$.

We now consider $\SPP_{\RR_+}(5)$. The check described in Subsection~\ref{subsec:small sumsets} would have found any additional sum-product pairs with $|A+A|\le 11$; there are none.
To finish, we need to consider the pairs $(i,j)$ with $12\le i \le 14$ and $9\le j \le 11$. As $j\le 11=3(5)-4$, the hypotheses of Theorem~\ref{thm:small product} are satisfied. Ergo, $A$ is a subset of a geometric progression with $7$ elements, say $Y=\{1,r,r^2,\dots,r^6\}$, with $1<r<2$ and irrational. Moreover, $r$ is a root of one of the $12$ polynomials
\begin{multline*}
  x^4-x^2-1,x^2-x-1,x^3-x-1,x^3-x^2-1,x^3-x^2-x-1,\\
  x^4-x-1,x^4-x^3-x^2-x-1,x^4+x^3-x^2-x-1,x^5-x-1,\\
  x^5-x^3-x^2-x-1,x^5-x^4-x^3-x^2-x-1,x^5+x^4-x^2-x-1.
\end{multline*}
Exhaustive computer search through the possible sets gives the following examples, and no other sum-product pair occurs.
\[
\begin{array}{cccc}
  |A+A| & |AA| & A & \text{minimal polynomial of $r$} \\ \hline
  13 & 9 & \left\{1,r,r^2,r^3,r^4\right\} & r^2-r-1 \\
  13 & 10 & \left\{1,r^2,r^3,r^4,r^5\right\} & r^2-r-1 \\
  13 & 11 & \left\{1,r,r^3,r^4,r^5\right\} & r^3-r^2-r-1 \\
  14 & 9 & \left\{1,r,r^2,r^3,r^4\right\} & r^3-r-1 \\
  14 & 10 & \left\{1,r,r^2,r^3,r^5\right\} & r^4+r^3-r^2-r-1 \\
  14 & 11 & \left\{1,r^2,r^3,r^4,r^6\right\} & r^5+r^4-r^2-r-1 \\
\end{array}
\]

We now consider $\SPP_{\RR_+}(6)$. The check described in Subsection~\ref{subsec:small sumsets} would have found any additional sum-product pairs with $|A+A|\le 14$; there are none.
To finish, we need to consider the pairs $(i,j)$ with $15\le i \le 20$ and $11\le j \le 14$. As $j\le 14=3(6)-4$, the hypotheses of Theorem~\ref{thm:small product} are satisfied. Ergo, $A$ is a subset of a geometric progression with $9$ elements, say $Y=\{1,r,r^2,\dots,r^8\}$, with $1<r<2$ and irrational. Moreover, $r$ is a root of one of $32$ polynomials.
Exhaustive computer search through the possible sets ($32$ values of $r$, and $\binom{9}{6}$ subsets of $Y$) gives the following examples, and no other sum-product pair occurs.
\[
\begin{array}{cccc}
  |A+A| & |AA| & A & \text{minimal polynomial of $r$} \\ \hline
  18 & 11 & \left\{1,r,r^2,r^3,r^4,r^5\right\} & r^2-r-1 \\
  18 & 12 & \left\{1,r^2,r^3,r^4,r^5,r^6\right\} & r^3-r-1 \\
  18 & 13 & \left\{1,r^3,r^4,r^5,r^6,r^7\right\} & r^3-r-1 \\
  18 & 14 & \left\{1,r^2,r^3,r^5,r^6,r^7\right\} & r^3-r-1 \\
  19 & 11 & \left\{1,r,r^2,r^3,r^4,r^5\right\} & r^3-r-1 \\
  19 & 12 & \left\{1,r,r^2,r^3,r^4,r^6\right\} & r^2-r-1 \\
  19 & 13 & \left\{1,r,r^3,r^4,r^5,r^6\right\} & r^4+r^3-r^2-r-1 \\
  19 & 14 & \left\{1,r^2,r^3,r^5,r^6,r^7\right\} & r^4+r^3-r^2-r-1 \\
  20 & 11 & \left\{1,r,r^2,r^3,r^4,r^5\right\} & r^4+r^3-r^2-r-1 \\
  20 & 12 & \left\{1,r,r^2,r^3,r^4,r^6\right\} & r^5+r^4-r^2-r-1 \\
  20 & 13 & \left\{1,r,r^2,r^3,r^4,r^7\right\} & r^6+r^5+r^4-r^3-r^2-r-1 \\
  20 & 14 & \left\{1,r,r^2,r^4,r^5,r^8\right\} & r^7+r^6+r^5-r^3-r^2-r-1 \\
\end{array}
\]
\end{proof}

\subsection{How Much Work Will \texorpdfstring{$\SPP(7)$}{SPP(7)} be?}
We have already performed the search that confirms that there are no additional sum-product pairs $(i,j)$ with $i\le 17$ in $\SPP(7)$, beyond those shown in Figure~\ref{fig:tiny n}. The SEZ confirms that there are no additional pairs in $\SPP(7)$ with $j\le 17$. While that sort of work was sufficient for $n\le 6$, for $n=7$ we have four additional points to consider:
    $(21,18),(18,20),(18,21),(19,20)$.
We have no examples that give these as sum-product pairs, nor any easy attack to show that they are not sum-product pairs.

\section{Problems Warranting Further Study}
\label{sec:questions}

Aside from the conjectures called out in the text above, there are several problem areas that have arisen while building the dataset and the visualizations.

\subsection{Types of Addition and Multiplication Tables}
The following considerations may be helpful in proving the exact value of $\SPP(n)$ for $n\ge 7$.

If one orders $A=\{a_1<\cdots<a_n\}$ and constructs the addition table, one readily observes that each row and each column is strictly increasing. In this way, each set $A$ of $n$ real numbers induces a weak ordering on the set of pairs $\{(i,j) : 1\le i \le j \le n\}$. Namely, $(i,j)=(k,\ell)$ if $a_{i}+a_{j}=a_{k}+a_{\ell}$ and $(i,j)<(k,\ell)$ if $a_i+a_j<a_k+a_\ell$.

With this motivation, we call a weak ordering of $\{(i,j) : 1\le i \le j \le n\}$ a \emph{prototype of order $n$} if $(i,j)<(i+1,j)$ and $(i,j)<(i,j+1)$ for all $i,j$ for which this makes sense. How many prototypes of order $n$ are there? We computed the values for $n\le 6$,
\[1, 1, 3, 39, 2\,905, 1\,538\,369\]
and Howroyd~\cite[OEIS sequence \seqnum{A376162}]{oeis} has extended the known values to $n=20$.

We call a prototype a \emph{type of addition table} if it actually arises from a set $A$ as described above. How many types of addition tables are there? The sequence starts $1, 1, 3, 25, 477, \ldots$, and is sequence $\seqnum{A378609}$ in the OEIS.

It is not hard to prove that the types of addition tables are exactly the types of multiplication tables. To show that $(19,20)\not\in \SPP(7)$, one can then consider the types of addition tables that lead to $|A+A|=19$, and it is a linear algebra problem for each type to express $a_1,\dots,a_n$ in terms of a few parameters. Then one considers the possible types of multiplication tables that lead to $|AA|=20$, and the type specifies that certain pairs of products are equal and that other pairs of products satisfy strict inequalities. This leads to a system of quadratic equations and inequations in the parameters that is typically easy to solve. The difficulty here is in enumerating the types, and then automating the solving of the resulting equations and inequalities. This author was able to get this approach working for $n=4$, but not for $n=5$, much less for $n=7$.

\subsection{Varying the Domain}
At the beginning of this project, the author chose to consider sets of $n\le 32$ positive integers almost arbitrarily. ``Positive integers'' seemed like the most popular domain, and he believed that $n=20$ would be large enough to begin seeing the sum-product force push points north and east, and decided to collect data to $n=32$ as a safety margin.

It would be instructive to find pairs $(x,y)$ that prove that $\SPP_{\ZZ}(n)$, $\SPP_{\NN_0}(n)$, $\SPP_{\RR}(n)$, etc., are distinct.

It seems unlikely that
 \[
 \SPP_{\RR_+}(n)=\SPP(n) \cup \left( [\tbinom{n+1}2-n+3,\tbinom{n+1}2-1] \times [2n-1,3n-4]\right), \]
for all $n$, but this holds at least for $n\le 6$. We believe that $(23,13)$ is not in $\SPP_{\RR_+}(7)$.

\subsection{The River of Ignorance}
Uniformly random subsets of $[N]$ tend to give normalized sum-product pairs along the north edge, while subsets of the set of divisors gives pairs along the eastern edge. There is a ``River of Ignorance'' running between the two that we have filled in only by the shotgun methods of considering all translates of the sets already in our dataset and by trying to adjoin many possible numbers to each set already in our dataset. We have no reliable way of generating sets with, for example $|A+A|\sim n^{1.8}$ and $|AA|\sim n^{1.9}$. Of course, if the sum-product conjecture is true, then there is no such way that works for large $n$.

More generally, none of our structured techniques of generating sets consistently gave sets with normalized sum-product pair near the segment connecting $(1,1.5),(2,2)$. What types of sets  populate the author's River of Ignorance?

The lower right of our images, just above the flagpole, is unnaturally unpopulated. This is likely because the author doesn't know how to efficiently generate many sets with $3n-3\le |AA|\le 4n$ that are not Sidon sets.

\subsection{Tiny \texorpdfstring{$n$}{n}}
Generate plausible conjectures for $\SPP(n)$ for larger $n$, say $n\le 32$. This will require a more systematic and extensive search than the one described in this work.

\subsection{Infrastructure}
Create a website that stores for each $(i,j,n)$ the extreme $n$-element sets with sumset of size $i$ and product set of size $j$. We would want to record both the smallest diameter example, and the smallest maximum (if different). We would also want $X=\NN$ of course, but maybe also $\NN_0$ and $\ZZ$. This should be plausible for $n\le 50$. Maybe up to $n\approx 100$ if one only stores the extremal $i,j$ for each $n$. There should be a user interface allowing people to submit additional sets, with the website checking if they add to the database. For each such set, we would want to store the submitter, the date of submission, and maybe a note as to how the set was generated/discovered.

\subsection{The True Size of \texorpdfstring{$\SPP(n)$}{SPP(n)}}
Is $\abs{\SPP(n)}\sim \frac14 n^4$? This is sequence \seqnum{A378623} in the OEIS. The only improvements in the upper bound are $O(n^3)$, like the SEZ and the Solymosi Void. The only lower bound that this author is aware of is mentioned above: for each $t$ between $2n-1$ and $n(n+1)/2$, there is a set $A$ of $n$ positive integers with $|A+A|=t$, and a set with $|AA|=t$.

The lower bound for the first terms of the sequence given by our dataset is given in the following table. The values for $n\le 6$ are proven exact, the others are merely lower bounds.
\[\begin{array}{|c|c|c|c|c|c|c|c|c|} \hline
 n & 1 & 2 & 3 & 4 & 5 & 6 & 7 & 8 \\
|\SPP(n)| & 1 & 1 & 3 & 10 & 27 & 69 & 153 & 305 \\ \hline
n & 9 & 10 & 11 & 12 & 13 & 14 & 15 & 16 \\
|\SPP(n)| & 543 & 914 & 1444 & 2185 & 3198 & 4520 & 6233 &
8400 \\ \hline
n & 17 & 18 & 19 & 20 & 21 & 22 & 23 & 24 \\
|\SPP(n)| & 11081 & 14360 & 18342 & 23065 & 28675 & 35269 &
42935 & 51705 \\ \hline
n & 25 & 26 & 27 & 28 & 29 & 30 & 31 & 32 \\
|\SPP(n)| & 61822 & 73306 & 86294 & 100835 & 117196 & 135540
& 155890 & 178547 \\ \hline
\end{array}
\]

\subsection{Real Small Product Sets}
What is the smallest possible value of $|A+A|$, where $A$ is a set of $n$ positive real numbers with $|AA|\le 3n-4$? Is it achieved by $\{\varphi^i : 1\le i \le n\}$, where $\varphi=(1+\sqrt5)/2$?

\subsection{Freiman's \texorpdfstring{$(3n-3)$}{(3n-3)}-Theorem}
Freiman has a theorem that describes the sets of integers with $|A+A|=3n-3$. This may allow one to describe the sum-product pairs $(i,j)\in\SPP(n)$ with $j=3n-3$. Less likely, but also plausible, is to describe those sum-product pairs in $\SPP_{\RR_+}(n)$.

\subsection{The Other Flagpole}
In Figure~\ref{fig:everythingK}, one can just barely see a flagpole structure along the left edge. Is there some truth along the lines of: if $|A+A|\le 3n-4$ and $|AA| \le f(|A+A|)$, then $A$ is an arithmetic progression? Of course, this is vacuously true if $f(x)=1$, but one hopes for a choice of $f$ that makes the hypothesis satisfiable for arbitrarily large sets.

What are the triples $(\alpha,\beta,\gamma)$ with $|A+A|^\alpha|AA|^\beta \gg |A|^\gamma$ for every finite set $A$ with $|A|\ge 3$? To account for flagpoles, one likely wants to exclude arithmetic and geometric progressions, and Sidon sets and multiplicative Sidon sets.

\section{In Defense of the Sum-Product Conjecture}
\label{sec:defending}
We do not see evidence in support of the Sum-Product Conjecture in our dataset. Nevertheless, we now consider how to defend the conjecture, and criticize the utility of our dataset.

First, the dataset only contains sets with $32$ or fewer elements. In additive combinatorics, $32$ is often large enough to begin seeing asymptotic phenomena, but in this case we are considering product sets. Product sets involve prime decompositions, and seeing the typical behavior of product sets may require one to consider sets large enough to involve typical prime numbers. In prime number theory, 32 is microscopically small.

The dataset is not complete. Perhaps a more complete dataset will have better examples with smaller sumset and product sets for $n<32$, and we will see a gradual increase in our upper bound on
  \[\alpha_n \coloneqq \min_{\substack{A\subseteq\NN \\ |A|=n}} \max\{|A+A|,|AA|\},\]
instead of the gradual decrease that is apparent today. That is, let $u_n$ be the upper bound on $\alpha_n$ that arises from our dataset; further examples can only reduce the value of each $u_n$. But the suggestion that the Sum-Product Conjecture is false comes not from the values of $u_n$ but from the trend $u_{28}>u_{29}>u_{30}>u_{31}>u_{32}$. This downward trend suggests that $\alpha_n$ is also trending downward, away from 2. However, further examples may produce new (lower) values of $u_n$ with a different trend, say $u_{28}<u_{20}<\cdots < u_{32}$, which would no longer impugn Erd\H{o}s's conjecture.

Sets of $y$-friable numbers\footnote{Often called ``smooth numbers''}, numbers all of whose prime factors are at most $y$, give many of the most extreme examples in the dataset. These are at least partially understood and give $\displaystyle \min_{|A|=n} \max\{|A+A|,|AA|\}\to 2$, but at a glacial pace. In particular, a pace too slow to see with $n\le 32$.

The Solymosi Void contains many normalized sum-product pairs, including points with $n=32$ along the lower envelope of discovered points. And yet we know that there are no limit points of $\NSPP([3,\infty))$ inside the Solymosi Void. It follows that $n=32$ is not large enough for the sum-product pressure to kick in.

\section*{Acknowledgements}
I thank the attendees at the New York Number Theory Seminar for many helpful suggestions for building out the dataset and for presenting it. Extended conversations with Moshe Newman were particularly helpful. The referees provided considerable help and eliminated a major error, for which I am exceedingly grateful.


\begin{thebibliography}{10}

  \bibitem{Bloom}
  Thomas~F. Bloom.
  \newblock Control and its applications in additive combinatorics, 2025.

  \bibitem{2003.Chang}
  Mei-Chu Chang.
  \newblock The {E}rd{\H o}s-{S}zemer\'edi problem on sum set and product set.
  \newblock {\em Ann. of Math. (2)}, 157(3):939--957, 2003.

  \bibitem{Rice}
  Ginny~Ray Clevenger, Haley Havard, Patch Heard, Andrew Lott, Alex Rice, and
  Brittany Wilson.
  \newblock The sum-product problem for small sets.
  \newblock {\em Involve}, 18(1):165--180, 2025.

  \bibitem{1976.Erdos}
  P.~Erd\H{o}s.
  \newblock Some recent problems and results in graph theory, combinatorics and
  number theory.
  \newblock In {\em Proceedings of the {S}eventh {S}outheastern {C}onference on
    {C}ombinatorics, {G}raph {T}heory, and {C}omputing ({L}ouisiana {S}tate
    {U}niv., {B}aton {R}ouge, {L}a., 1976)}, volume No. XVII of {\em Congress.
    Numer.}, pages 3--14. Utilitas Math., Winnipeg, MB, 1976.

  \bibitem{ErdosSzemeredi}
  P.~Erd\H{o}s and E.~Szemer\'edi.
  \newblock On sums and products of integers.
  \newblock In {\em Studies in pure mathematics}, pages 213--218. Birkh\"auser,
  Basel, 1983.

  \bibitem{Freiman2}
  Gregory Freiman, Marcel Herzog, Patrizia Longobardi, and Mercede Maj.
  \newblock Small doubling in ordered groups.
  \newblock {\em J. Aust. Math. Soc.}, 96(3):316--325, 2014.

  \bibitem{Freiman}
  G.~A. Fre\u~iman.
  \newblock The addition of finite sets. {I}.
  \newblock {\em Izv. Vys\v s. U\v cebn. Zaved. Matematika},
  1959(6(13)):202--213, 1959.

  \bibitem{UPINT}
  Richard~K. Guy.
  \newblock {\em Unsolved problems in number theory}.
  \newblock Problem Books in Mathematics. Springer-Verlag, New York, third
  edition, 2004.

  \bibitem{Nathansonsumsetsizes}
  Melvyn~B. Nathanson.
  \newblock Problems in additive number theory, {V}{I}: {S}izes of sumsets, 2024.

  \bibitem{oeis}
  {OEIS Foundation Inc.}
  \newblock The {O}n-{L}ine {E}ncyclopedia of {I}nteger {S}equences, 2024.
  \newblock Published electronically at \url{http://oeis.org}.

  \bibitem{Solymosi}
  J\'ozsef Solymosi.
  \newblock Bounding multiplicative energy by the sumset.
  \newblock {\em Adv. Math.}, 222(2):402--408, 2009.

\end{thebibliography}

\end{document}